\documentclass{article}
\usepackage{amsmath, amsfonts, amssymb, amsthm, bbm}

\newtheorem*{thma}{Theorem}
\newcommand{\reff}[1]{(\ref{#1})}

\newcommand {\cL} {\mathcal{L}}

\newcommand {\cX} {\mathcal{X}}
\newcommand {\cY} {\mathcal{Y}}

\newcommand{\vs} {\varsigma}
\begin{document}
\title{On Seneta's Constants for the Supercritical Bellman-Harris Process with $E(Z_+ \log Z_+) = \infty$}
\author{Wolfgang P. Angerer\\Goethe-Universit\"at\\Frankfurt am Main\\Germany}
\date{}
\maketitle
\newpage
\begin{abstract}
\noindent For a finite mean supercriticial Bellman-Harris process, let $Z_t$ be the number of particles at time $t$. There exist numbers $\chi_t$ (the Seneta constants) such that $\chi_t Z_t$ converges almost surely to a non-degenerate limit. Furthermore, $\chi_t \propto e^{-\beta t} \cL(e^{-\beta t})$, where $\beta$ is the Malthusian parameter, and $\cL$ is slowly varying at zero. We obtain a characterisation of the slowly varying part of the Seneta constants under the assumption that the life-time distribution of particles is strongly non-lattice.\\

\noindent Keywords: {\sc Branching Processes; Renewal Theory} \\

\noindent Mathematics Subject Classification (2000): 60J80
\end{abstract}

\newpage
\section{Introduction}
We consider a supercritical Bellman-Harris process $\{Z_t\}_{t\geq 0}$ with offspring distribution $\{\pi_k\}_{k = 0}^{\infty}$ and life-time distribution $G$. In words, $G(t)$ is the probability that a newborn particle survives at least until time $t$, and $\pi_k$ is the probability that once it splits into a number $Z_+$ of progeny, it will split into exactly $k$ of these. We denote by $f(s)$ the corresponding generating function (PGF) $f(s) := E(e^{sZ_+}) = \sum_{k = 0}^{\infty} \pi_k s^k$, and set
\begin{equation}
h(s) := \frac{1 - f(s)}{1 - s} \:.\label{h}
\end{equation}
Let $\mu := h(1)$ be finite. It is known that there exist \lq constants' $\chi_t$ (the Seneta constants) such that, on the set of non-extinction, $\chi_t Z_t$ converges almost surely to a non-degenerate random variable $Z$. With $F_t(s)$ the PGF of the distribution of particle numbers at time $t$, and $F_{-t}(s)$ its inverse, an immediate candidate for the $\chi_t$'s is
\begin{equation}
\chi_t :=: \chi_t(\vs) := - \log F_{- t}(\vs)
\end{equation}
for some $\vs \in (q, 1)$, where $q$ is such that $h(q) = 1$ (such a $q$ exists and is unique because of supercriticality). This is because the Laplace transform of the random variable $Z$ (with $y$ as the dummy variable) is 
\begin{equation}
R(y) :=: R_{\vs}(y) := \lim_{t\to\infty} F_t(e^{- \chi_t(\vs)y})\:, \label{ls}
\end{equation}
which by definition of $\chi_t$ equals $\vs$ for $y = 1$. Hence $Z$ is non-trivial in the sense that its Laplace transform is neither 0 nor 1. To get a feeling for how quickly the $\chi_t$'s tend to zero, recall that $Z_t$ grows essentially as $e^{\beta t}$ as $t\to\infty$, where $\beta$ is the Malthusian parameter, that is,
\begin{displaymath}
\int_0^{\infty} e^{-\beta t} \, dG(t) = \frac{1}{\mu} \:.
\end{displaymath}
It is therefore natural to conjecture that $\chi_t = e^{-\beta t} \cL(e^{-\beta t})$ for some slowly varying (in $e^{-\beta t}$) function $\cL$. Under this assumption, it is easy to derive an equation for the Laplace transform $R$ of the random variable $Z$: Since $F_t$ fulfills the integral equation \cite{an}
\[
F_t(s) = \big(1 - G(t)\big)s + \int_0^t f \circ F_{t-u}(s)\,dG(u)\:,
\]
it follows immediately that
\begin{eqnarray}
R(y)&=&\lim_{t\to\infty} \big(1- G(t)\big) e^{-\chi_t y} + \int_0^t f \circ F_{t-u}(e^{-\chi_{t-u}(\chi_t/\chi_{t-u}) y})\,dG(t)\nonumber\\
&=&\int_0^{\infty} f \circ R(y e^{-\beta t})\,dG(t)\:,\label{leq}
\end{eqnarray}
because of dominated convergence and the fact that $\chi_t/\chi_{t-u} = e^{-\beta u}$ except for a factor which tends to 1 as $t\to\infty$. Furthermore, since $\chi_t e^{-\beta t}$ is certainly very close to zero for $t$ large enough, we may approximate $1-R(e^{-\beta t})\sim 1 - F_t(e^{-\chi_t e^{-\beta t}})$ by $E(Z_t)\chi_t e^{-\beta t}\sim \chi_t$, so that $1-R(e^{-\beta t})$ would also be a natural guess at the value of the Seneta constants. (By $\sim$ we mean that the ratio of the two quantities is bounded from above and away from zero.) This guess is indeed a good one, as has been established by Schuh \cite{sch} for the case of $G$ being non-lattice, and makes it natural for us to consider
\begin{equation}
\cX(t) := e^{\beta t}\big(1 - R(e^{-\beta t})\big) = \int_0^{\infty} h \circ R(e^{-\beta (t+u)})\cX(t+u)e^{-\beta u}\,dG(u)\:. \label{tvs}
\end{equation}
Our idea is to assume $h$, $R$, and $G$ as given and treat Equation \reff{tvs} as an equation in the unknown function $\cX$. It will turn out that under the assumption that $G$ is {\it strongly non-lattice}, which is to say that
\[
\liminf_{\vert\theta\vert\to\infty}\left\vert 1 - \int_0^{\infty} e^{\sqrt{-1}\theta t}\,dG(t)\right\vert > 0\:,
\]
we will need no more than elementary renewal theory to derive the following
\begin{thma} \label{main} Suppose the life-time distribution $G$ of particles is strongly non-lattice, and that $E(Z_+ \log Z_+) = \infty$. Then, with
\begin{equation}
\nu^{-1} = \mu\int_0^{\infty} t e^{-\beta t}\,dG(t)\:,\label{nu}
\end{equation}
we have
\begin{equation}
\log \cX(t) \simeq \frac{\nu}{\mu}\,\int_0^t \big(\mu - h(1-e^{-\beta u})\big)\,du\:,\label{gist}
\end{equation}
where $\simeq$ means that the ratio of both sides tends to $1$ as $t\to\infty$.
\end{thma}
\noindent The gist of the theorem is that both sides of \reff{gist} are asymptotically equivalent with a constant of proportionality equal to $1$ instead of only $O(1)$. Note that $\cX$ does not depend anymore on $\vs$ on this scale. We also observe that under what might be called Uchiyama's \cite{u} condition:
\[
\mu - h(1-s) = (-\log s)^{-\alpha} \cL( -\log s)\:,
\]
where $\alpha\geq 0$ and $\cL$ is slowly varying at infinity, the theorem reduces to
\[
\log \cX(t) \simeq \frac{\nu}{\mu}\,\int_0^t (\beta u)^{-\alpha} \cL(\beta u)\,du \simeq \frac{\nu}{(1-\alpha)\beta^{\alpha}\mu}\, t^{1-\alpha} \cL(\beta t),
\]
which, if we set $G(t) = \delta_{\tau}(t)$ (the Dirac mass at some splitting time $\tau$), $\beta = \log\mu/\tau$, and $\nu = \tau^{-1}$, is what Uchiyama's Theorem \cite{u} claims for the ordinary Galton-Watson process. (Time should then be measured in units of $\tau$.) Let us now turn to the
\section{Proof of the Theorem}
To begin, write
\begin{equation}
\cX(t) = \frac{1}{\mu}\int_{u=t}^{\infty} h\circ R(e^{-\beta u})\cX(u)\,dG_{\beta}(u-t)\:,
\end{equation}
with
\begin{equation}
G_{\beta}(t) := \mu \int_0^t e^{- \beta u}\,dG(u)\:,\label{gb}
\end{equation}
and introduce the \lq renewal function\rq\ 
\begin{equation}
U_{\beta}(t) := \sum_{i = 1}^{\infty} G_{\beta}^{\ast i}(t)\:,
\end{equation}
where $G_{\beta}^{\ast i}$ is for the $i$-fold convolution of $G_{\beta}$ with itself. (The standard definition would be to also include a Dirac mass at zero.) By renewal theory \cite{f}, we have for $t>0$,
\begin{equation}
1 + U_{\beta}(t) = \nu t + \tilde{U}_{\beta}(t)\:,\label{rnt}
\end{equation}
where $\tilde{U}_{\beta}(t)\geq 0$, $\nu$ is as given in the theorem, and
\begin{equation}
\int_0^{\infty} d\tilde{U}_{\beta}(t) = \frac{\nu^2}{2}\,\int_0^{\infty} t^2\,dG_{\beta}(t)\:.
\end{equation}
Since $G_{\beta}$ has an exponential tail, and is strongly non-lattice together with $G$, $\tilde{U}_{\beta}$ has an exponential (right) tail as well \cite{st}. Consider now
\begin{eqnarray}
\lefteqn{\int_{u=0}^t \cX(u)\,dU_{\beta}(u) = \frac{1}{\mu}\int_{u=0}^t \int_{v=u}^{\infty} h\circ R(e^{-\beta v})\cX(v)\,dG_{\beta}(v-u)\,dU_{\beta}(u)}\nonumber\\
&=& \frac{1}{\mu}\int_{v=0}^{\infty} \int_{u=0}^{v \wedge t} h\circ R(e^{-\beta v})\cX(v)\,dG_{\beta}(v-u)\,dU_{\beta}(u)\nonumber\\
&=& \frac{1}{\mu}\int_{v=0}^t \int_{u=0}^v h\circ R(e^{-\beta v})\cX(v)\,dG_{\beta}(v-u)\,dU_{\beta}(u)\nonumber\\
&&+\:\frac{1}{\mu}\int_{v=t}^{\infty} \int_{u=0}^t h\circ R(e^{-\beta v})\cX(v)\,dG_{\beta}(v-u)\,dU_{\beta}(u)\nonumber\\
&=& \frac{1}{\mu}\int_{v=0}^t h\circ R(e^{-\beta v})\cX(v)\,d(G_{\beta}\ast U_{\beta})(v)\nonumber\\
&&+\:\frac{1}{\mu}\int_{v=t}^{\infty} h\circ R(e^{-\beta v})\cX(v)\,d\big(G_{\beta}\ast U_{\beta}(v) - G_{\beta}\ast U_{\beta}(v-t)\big)\nonumber\\
&&-\:\frac{1}{\mu}\int_{v=t}^{\infty} \int_{u=t}^{v} h\circ R(e^{-\beta v})\cX(v)\,dG_{\beta}(v-u)\,d\big(\tilde{U}_{\beta}(u) - \tilde{U}_{\beta}(u-t)\big)\:,\nonumber
\end{eqnarray}
since, by Equation \reff{rnt},
\[
U_{\beta}(u+t) = \nu t + U_{\beta}(u) + \tilde{U}_{\beta}(u+t) - \tilde{U}_{\beta}(u)\:.
\]
But
\[
G_{\beta}\ast U_{\beta}(t) = U_{\beta}\ast G_{\beta}(t) = U_{\beta}(t) - G_{\beta}(t)
\]
by definition of the renewal function, so we obtain
\begin{eqnarray}
\lefteqn{\int_{u=0}^t \cX(u)\,dU_{\beta}(u) = \frac{1}{\mu}\int_{v=0}^t h\circ R(e^{-\beta v})\cX(v)\,dU_{\beta}(v)}\nonumber\\
&&-\: \frac{1}{\mu}\int_{v=0}^t h\circ R(e^{-\beta v})\cX(v)\,dG_{\beta}(v)\nonumber\\
&&+\:\frac{1}{\mu}\int_{v=t}^{\infty} h\circ R(e^{-\beta v})\cX(v)\,d\big(\tilde{U}_{\beta}(v) - \tilde{U}_{\beta}(v-t)\big)\nonumber\\
&&-\:\frac{1}{\mu}\int_{v=t}^{\infty} h\circ R(e^{-\beta v})\cX(v)\,dG_{\beta}(v)\:+\:\frac{1}{\mu}\int_{v=t}^{\infty} h\circ R(e^{-\beta v})\cX(v)\,dG_{\beta}(v-t)\nonumber\\
&&-\:\frac{1}{\mu}\int_{v=t}^{\infty} \int_{u=t}^v h\circ R(e^{-\beta v})\cX(v)\,dG_{\beta}(v-u)\,d\big(\tilde{U}_{\beta}(u) - \tilde{U}_{\beta}(u-t)\big)\:.\nonumber
\end{eqnarray}
The second and forth term on the right-hand side of this equation add to $-\cX(0)$, by Equation \reff{rnt}. The fifth term, by the same equation, is simply $\cX(t)$. As for the remaining term, we have
\begin{eqnarray}
\lefteqn{\frac{1}{\mu}\int_{v=t}^{\infty} \int_{u=t}^v h\circ R(e^{-\beta v})\cX(v)\,dG_{\beta}(v-u)\,d\big(\tilde{U}_{\beta}(u) - \tilde{U}_{\beta}(u-t)\big)}\nonumber\\
&=&\int_{u=t}^{\infty}\frac{1}{\mu}\int_{v=u}^{\infty} h\circ R(e^{-\beta v})\cX(v)\,dG_{\beta}(v-u)\,d\big(\tilde{U}_{\beta}(u) - \tilde{U}_{\beta}(u-t)\big)\nonumber\\
&=&\int_{u=t}^{\infty}\cX(u)\,d\big(\tilde{U}_{\beta}(u) - \tilde{U}_{\beta}(u-t)\big)\:,\nonumber
\end{eqnarray}
so that
\begin{eqnarray}
\lefteqn{\int_{u=0}^t \cX(u)\,dU_{\beta}(u) = \frac{1}{\mu}\int_{u=0}^t h\circ R(e^{-\beta u})\cX(u)\,dU_{\beta}(u)}\nonumber\\
&&+\:\cX(t) - \cX(0) -\int_{u=t}^{\infty} \frac{\mu - h\circ R(e^{-\beta u})}{\mu}\,\cX(u)\,d\big(\tilde{U}_{\beta}(u) - \tilde{U}_{\beta}(u-t)\big)\:.\nonumber
\end{eqnarray}
If now we make use of \reff{rnt} once more, we finally obtain
\begin{equation}
\frac{\nu}{\mu}\int_{u=0}^t \big(\mu - h\circ R(e^{-\beta u})\big)\cX(u)\,du = \cX(t) - \cX(0) + \tilde{\cX}(t) - \tilde{\cX}(0)\:,
\end{equation}
where
\begin{equation}
\tilde{\cX}(t) := \int_{u=t}^{\infty} \frac{\mu - h\circ R(e^{-\beta u})}{\mu}\,\cX(u)\,d\tilde{U}_{\beta}(u-t)\:.
\end{equation}
With
\begin{equation}
\sigma(t) := \frac{\tilde{\cX}(t)}{\cX(t)} = \int_{u=0}^{\infty} \frac{\mu - h\circ R(e^{-\beta (t+u)})}{\mu}\,\frac{\cX(t+u)}{\cX(t)}\, d\tilde{U}_{\beta}(u)\geq 0\:,\label{sigm}
\end{equation}
we can write down a \lq solution\rq\ for this as
\begin{equation}
\cX(t) = (1-\vs)\,\frac{1+\sigma(0)}{1+\sigma(t)} \exp\left(\frac{\nu}{\mu}\,\int_0^t \frac{\mu - h\circ R(e^{-\beta u})}{1+\sigma(u)}\,du\right)\:.\label{form}
\end{equation}
We can assume that under $E(Z_+ \log Z_+)=\infty$, the integral in the exponent diverges (see, for instance, Corollary 3.1 in \cite{sch}, or Lemma 2 in \cite{u}). Therefore,
\[
\log \cX(t) \simeq \frac{\nu}{\mu}\,\int_0^t \frac{\mu - h\circ R(e^{-\beta u})}{1+\sigma(u)}\,du
\]
certainly if $\sigma$ is bounded for all sufficiently large $t$. We show that this is the case: Since $\cX(t)$ is slowly varying in $e^{-\beta t}$, it follows by Potter's Theorem \cite{bgt} that for arbitrary $A>1$, $\delta>0$, there exists $\tau_1 := \tau_1(A,\delta)$ such that
\[
\frac{\cX(t+u)}{\cX(t)}\leq A e^{\delta\beta u}
\]
for all $t\geq\tau_1$. Thus we obtain for all such $t$, from Equation \reff{sigm},
\begin{equation}
\sigma(t) \leq A\,\frac{\mu - h\circ R(e^{-\beta t})}{\mu}\, \int_{u=0}^{\infty} e^{\delta\beta u}\,d\tilde{U}_{\beta}(u)\:,
\label{bsigm}\end{equation}
which is of order $\sim \mu - h\circ R(e^{-\beta t})$ for suitably chosen $\delta$, because $\tilde{U}_{\beta}$ has an exponential tail ($\delta < 1$ already suffices). In particular, $\sigma(t)\to 0$ as $t\to\infty$, which readily implies that even
\begin{equation}
\log \cX(t)\simeq \frac{\nu}{\mu}\,\int_0^t \big(\mu - h\circ R(e^{-\beta u})\big)\,du =:\cY(t)\:.\label{defy}
\end{equation}
Now, because of \reff{bsigm}, and due to the divergence of the integral in \reff{form}, there exists $\tau_2$ such that $\cX(t)\geq 1$ for all $t\geq \tau_2$, and
\[
\cY(t)-\cY(\tau_2) > \frac{\nu}{\mu}\, \int_{\tau_2}^t \big(\mu - h\circ(1 - e^{-\beta u})\big)\,du\:,
\]
by monotonicity of $h$. But $\cY(\tau_2)$ is finite, hence
\begin{equation}
\liminf_{t\to\infty} \frac{\cY(t)}{\nu\mu^{-1}\int_0^t \big(\mu - h\circ(1 - e^{-\beta u})\big)\,du}\geq 1\:.\label{lb}
\end{equation}
On the other hand, because $\cX(t)$ is slowly varying in $e^{-\beta t}$, we can choose $A>1$, $\delta > 0$, and $\tau_3<\infty$ such that $A e^{-\delta\beta\tau_3} \cX(\tau_3)<1$, and
\[
\cY(t)-\cY(\tau_3)<\frac{\nu}{\mu}\,\int_{\tau_3}^{t} \Big(\mu - h\circ \big(1 - e^{-\beta (1-\delta) u}\big)\Big)\,du
\]
for all $t \geq \tau_3$, by monotonicity of $h$, and Potter's Theorem again. But $\cY(\tau_3)$ is finite, hence
\begin{equation}
\limsup_{t\to\infty} \frac{\cY(t)}{\nu\mu^{-1}\int_0^t \big(\mu - h\circ(1 - e^{-\beta u})\big)\,du}\leq \frac{1}{1-\delta}\:,\label{ub}
\end{equation}
which together with \reff{lb} and the definition \reff{defy} of $\cY$ concludes the proof of the theorem, since $\delta$ is arbitrary. \hfill$\Box$
\medskip

{\bf Acknowledgment.} The author is happy to acknowledge financial support from DFG and NWO as part of a Dutch-German research project on random spatial models from physics and biology. He thanks Martin Hutzenthaler, Anton Wakolbinger, and everyone at ISMI for ongoing discussions.

\end{document}